# Hopf monoidal comonads

Dimitri Chikhladze,[1]   Stephen Lack,[2]  and Ross Street[1]

**1** Macquarie University

**2** University of Western Sydney

ABSTRACT.  Alain Bruguières, in his talk [1], announced his work [2] with Alexis Virelizier and the second author which dealt with lifting closed structure on a monoidal category to the category of Eilenberg-Moore algebras for an opmonoidal monad. Our purpose here is to generalize that work to the context internal to an autonomous monoidal bicategory. The result then applies to quantum categories and bialgebroids.

## 1. Introduction

Early workers on monoidal categories recognized the importance of functors $F : \mathcal{V} \to \mathcal{W}$ not preserving the tensor, rather, having only coherent natural families of morphisms

$$\phi : F\,X \otimes F\,Y \longrightarrow F(X \otimes Y)\,.$$

There is also a morphism $\phi_0 : I \to F\,I$ connecting the tensor units. We follow Eilenberg and Kelly [3] in using the term *monoidal functor* for this structure on $F$. A monoidal category $\mathcal{V}$ is *right closed* when each pair of objects $X$ and $Z$ has a right internal hom $Z^X$; that is, there is an evaluation morphism $\mathrm{ev} : X \otimes Z^X \to Z$ determining a bijection

$$\mathcal{V}\big(Y, Z^X\big) \cong \mathcal{V}(X \otimes Y, Z)\,, \quad f \mapsto \mathrm{ev}\,(1_X \otimes f)\,.$$

The monoidal structure $\phi$ on $F$ can be expressed in terms of right closed structure

$$\phi^r : F\big(Z^X\big) \longrightarrow (F\,Z)^{F\,X}$$

defined by the equality

$$F\,X \otimes F\big(Z^X\big) \xrightarrow{1 \otimes \phi^r} F\,X \otimes (F\,Z)^{F\,X} \xrightarrow{\mathrm{ev}} F\,Z$$
$$= \ F\,X \otimes F\big(Z^X\big) \xrightarrow{\phi} F\big(X \otimes Z^X\big) \xrightarrow{F(\mathrm{ev})} F\,Z\,.$$

It was also recognized that tensor preserving (up to isomorphism) was an important special case; these are the *strong monoidal functors* where "strong" means that $\phi$ and $\phi_0$ are invertible. Internal hom preserving was considered too rare to warrant terminology. However, it does occur in interesting cases such as the forgetful functor from $G$-sets to sets where $G$ is a group, a fact that distinguishes groups $G$ from monoids. This observation was important to both the papers [10] and [19] in defining "Hopf concepts" (in the sense of existence of antipodes). Following [12], we call a monoidal functor *strong right closed* when $\phi^r$ and $\phi_0$ are invertible. While strong monoidal functors do preserve duals, they are generally not strong right closed.

It was also recognized that tensor preserving (up to isomorphism) was an important special case; these are the *strong monoidal functors* where "strong" means that $\phi$ and $\phi_0$ are invertible. Internal hom preserving was considered too rare to warrant terminology. However, it does occur in interesting cases such as the forgetful functor from $G$-sets to sets where $G$ is a group, a fact that distinguishes groups $G$ from monoids. This observation was important to both the papers [10] and [19] in defining "Hopf concepts" (in the sense of existence of antipodes). Following [12], we call a monoidal functor *strong right closed* when $\phi^r$ and $\phi_0$ are invertible. While strong monoidal functors do preserve duals, they are generally not strong right closed.

Bruguières and Virelizier [4] began a program to generalize the theory of Hopf monoids in a braided monoidal category to appropriate monads on a not-necessarily-braided monoidal category. That paper looks at the case where the monoidal category is autonomous (duals exist). They define the concept of (right) antipode for an opmonoidal monad on a (right) autonomous monoidal category; such an opmonoidal monad with an antipode they call a Hopf monad. In later work [1-2] with the second author of the present paper they found a characterization of the existence of an antipode on the opmonoidal monad purely in terms of the monoidal structure of the autonomous monoidal category. (This corresponds to the fact for bimonoids that an antipode exists if and only if the fusion map is invertible.) The antipode-free formulation allowed them to extend the Hopf concept to opmonoidal monads on closed monoidal categories. For a right closed monoidal category they then showed that this Hopf condition on an opmonoidal monad was equivalent to the right closedness of the monoidal category of Eilenberg-Moore algebras together with the strong right closedness of the forgetful functor.

Day and the third author [12] defined a concept of quantum category in terms of monoidal comonads. They then define a quantum groupoid to be a quantum category equipped in a certain way with star-autonomous structure. Pastro and the third author [5] characterized when star-autonomous structure could be lifted to monoidal categories of Eilenberg-Moore coalgebras for a monoidal comonad.

The work of [1] immediately suggested the possibility of strengthening the quantum groupoid definition to a property of quantum category rather than extra structure. To prepare for this, in the present paper, we generalize the work of [1] to the level of pseudo-monoids in a monoidal bicategory. To suit our purpose we work with monoidal comonads on the pseudomonoids rather than opmonoidal monoids. The effect is that we need to take the concept of "coclosed" seriously since it is this structure that we wish to lift to the Eilenberg-Moore coalgebra pseudomonoid of the comonad.

We use the term "monoidale" rather than "pseudomonoid" or "monoidal object". After two sections establishing the context, in Section 4 we study the Hopf condition on monoidal comonads on monoidales in an autonomous monoidal bicategory. The revised definition of quantum groupoid in a reasonably general braided monoidal category appears in Section 5.

## 2. The bicategory of comonads

Let $\mathcal{M}$ denote a bicategory. We often ignore the associativity and unitivity constraints, writing as if $\mathcal{M}$ were a 2-category.

Essentially following [6] and [7], we consider a bicategory $\mathrm{Comnd}\mathcal{M}$ whose objects are pairs $(A, g)$ where $g = (g : A \to A, \delta : g \Rightarrow g\,g, \varepsilon : g \Rightarrow 1_A)$ is a comonad on the object $A$ of $\mathcal{M}$. A morphism $(k, \kappa) : (A, g) \to (A', g')$ in $\mathrm{Comnd}\mathcal{M}$ consists of a morphism $k : A \to A'$ and 2-cell $\kappa : k\,g \Rightarrow g'\,k$ in $\mathcal{M}$ satisfying the two conditions

$$\left( k\,g \xrightarrow{\kappa} g'\,k \xrightarrow{\delta\,k} g'\,g'\,k \right) = \left( k\,g \xrightarrow{k\,\delta} k\,g\,g \xrightarrow{\kappa\,g} g'\,k\,g \xrightarrow{g'\,\kappa} g'\,g'\,k \right),$$

$$\left( k\,g \xrightarrow{\kappa} g'\,k \xrightarrow{\varepsilon\,k} k \right) = \left( k\,g \xrightarrow{k\,\varepsilon} k \right).$$

A 2-cell $\sigma : (k, \kappa) \Rightarrow (l, \lambda) : (A, g) \to (A', g')$ in $\mathrm{Comnd}\mathcal{M}$ is a 2-cell $\sigma : k \Rightarrow l : A \to A'$ satisfying the condition

$$\left( k\,g \xrightarrow{\kappa} g'\,k \xrightarrow{g'\,\sigma} g'\,l \right) = \left( k\,g \xrightarrow{\sigma\,g} l\,g \xrightarrow{\lambda} g'\,l \right).$$

Here is an example of "doctrinal adjunction"[8] based on Theorem 1 of [9].

**Proposition 2.1.** *A morphism* $(k, \kappa) : (A, g) \to (A', g')$ *in* $\mathrm{Comnd}\mathcal{M}$ *has a left adjoint if and only if* $k$ *has a left adjoint and* $\kappa$ *is invertible in* $\mathcal{M}$.

*Proof.* Pseudofunctors, in particular the forgetful $\mathrm{Comnd}\mathcal{M} \to \mathcal{M}$, preserve adjunction so $k$ has a left adjoint if $(k, \kappa)$ does. The structure making this left adjoint into a morphism in $\mathrm{Comnd}\mathcal{M}$ is the mate of the inverse of $\kappa$. $\square$

Recall [7] that $\mathcal{M}$ admits the Eilenberg-Moore construction for comonads when the pseudofunctor $\mathcal{M} \to \mathrm{Comnd}\mathcal{M}$, taking $X$ to $(X, 1_X)$, has a right biadjoint $\mathrm{Comnd}\mathcal{M} \to \mathcal{M}$ whose value at $(A, g)$ is denoted by $A^g$. There is a pseudonatural equivalence

$$\mathcal{M}(X, A^g) \simeq \mathcal{M}(X, A)^{\mathcal{M}(X, g)}$$

where the right side is the category of Eilenberg-Moore coalgebras for the comonad $\mathcal{M}(X, g)$ on the category $\mathcal{M}(X, A)$; equally, the right side is the category $(\mathrm{Comnd}\mathcal{M})((X, 1_X), (A, g))$ whose objects we call $g$-*coalgebras*. Taking $X = A^g$ and evaluating the pseudonatural equivalence at the identity, we obtain the *universal* $g$-coalgebra $(u, \gamma) : (A^g, 1_{A^g}) \to (A, g)$ and remind the reader that this $u : A^g \to A$ has a right adjoint $u^*$ in $\mathcal{M}$.

**Corollary 2.2.** *If the morphism* $(k, \kappa) : (A, g) \to (A', g')$ *has a left adjoint in* $\mathrm{Comnd}\mathcal{M}$ *then the morphism* $k^\kappa : A^g \to A'^{g'}$ *has a left adjoint in* $\mathcal{M}$, *assuming the Eilenberg-Moore constructions exist.*

## 3. The bicategory of monoidales

Let $\mathcal{M}$ denote a monoidal bicategory [10]. We often ignore the various constraint equivalences, writing as if $\mathcal{M}$ were a Gray monoid; this is justified by [11].

A *monoidale* $E$ in $\mathcal{M}$ is what has previously been called a "pseudomonoid" or "monoidal object" (for example, see [12] and [13]). It consists of an object $E$, morphisms $p : E \otimes E \to E$ and $j : I \to E$, and invertible 2-cells $\alpha : p(p \otimes 1_E) \Rightarrow p(1_E \otimes p)$, $\lambda : p(j \otimes 1_E) \Rightarrow 1_E$ and $\rho : p(1_E \otimes j) \Rightarrow 1_E$ satisfying two axioms. For example, a monoidale in the cartesian monoidal bicategory Cat is a monoidal category.

For monoidales $E$ and $F$, a *monoidal morphism* $f = (f, \phi, \phi_0) : E \to F$ consists of a morphism $f : E \to F$ and 2-cells $\phi : p(f \otimes f) \Rightarrow f\, p$, $\phi_0 : j \Rightarrow f\, j$ in $\mathcal{M}$ satisfying three axioms. For example, a monoidal morphism $f : 1 \to F$ in Cat is a monoid in the monoidal category $F$.

The composite of monoidal morphisms $f = (f, \phi, \phi_0) : E \to F$ and $g = (g, \phi, \phi_0) : F \to G$ is defined to be $g\, f = (g\, f, \phi, \phi_0) : E \to F$ where

$$\phi = \left( p(g\, f \otimes g\, f) \cong p(g \otimes g)\,(f \otimes f) \xRightarrow{\phi(f \otimes f)} g\, p\,(f \otimes f) \xRightarrow{g\,\phi} g\, f\, p \right)$$

$$\text{and} \quad \phi_0 = \left( j \xRightarrow{\phi_0} g\, j \xRightarrow{g\,\phi_0} g\, f\, j \right).$$

A monoidal 2-cell $\sigma : f \Rightarrow f' : E \to F$ is a 2-cell in $\mathcal{M}$ satisfying $(\sigma\, p)\,\phi = \phi(p(\sigma \otimes \sigma))$ and $(\sigma\, j)\,\phi_0 = \phi_0$.

With obvious 2-cell compositions, this defines a bicategory $\mathrm{Mon}\mathcal{M}$ of monoidales in $\mathcal{M}$.

A monoidale in $\mathcal{M}^{\mathrm{co}}$ is a monoidale in $\mathcal{M}$ since the 2-cells $\alpha$, $\lambda$ and $\rho$ are replaceable by their inverses. However, an *opmonoidal morphism* $w = (w, \psi, \psi_0) : E \to F$ in $\mathcal{M}$ is a monoidal morphism in $\mathcal{M}^{\mathrm{co}}$; so we have a morphism $w : E \to F$ and 2-cells $\psi : w\, p \Rightarrow p(w \otimes w)$, $\psi_0 : w\, j \Rightarrow j$ in $\mathcal{M}$ satisfying three axioms.

A monoidal morphism $f = (f, \phi, \phi_0) : E \to F$ is called *strong* when $\phi_2$ and $\phi_0$ are invertible. Notice that a strong monoidal morphism yields an opmorphism by taking $\psi = \phi^{-1}$ and $\psi_0 = \phi_0{}^{-1}$. For a morphism $w$ with a right adjoint $w^*$, mateship provides a bijection between monoidal structures on $w^*$ and opmonoidal structures on $w$. Here is another example of "doctrinal adjunction" [8].

**Proposition 3.1.** *A monoidal morphism* $f = (f, \phi, \phi_0) : E \to F$ *has a right adjoint in* $\mathrm{Mon}\mathcal{M}$ *if and only if it is strong and* $f : E \to F$ *has a right adjoint in* $\mathcal{M}$. $\square$

Suppose $\mathcal{M}$ is right autonomous: for each object $A$ there is an object $A^{\circ}$ and an equivalence

$$\mathcal{M}(A \otimes B, C) \simeq \mathcal{M}(B, A^{\circ} \otimes C), \quad h \mapsto \hat{h}$$

pseudonatural in $B$ and $C$. It follows that the object assignment $A \mapsto A^{\circ}$ becomes a pseudofunctor $(-)^{\circ} : \mathcal{M}^{\mathrm{op}} \to \mathcal{M}$ and there are morphisms $e : A \otimes A^{\circ} \to I$ and $n : I \to A^{\circ} \otimes A$, extraordinary pseudonatural, and so on (for example, see [10]). It also follows that we have a pseudonatural equivalence pseudonatural in $B$ and $C$. It follows that the object assignment $A \mapsto A^\circ$ becomes a pseudofunctor $(-)^\circ : \mathcal{M}^{\mathrm{op}} \to \mathcal{M}$ and there are morphisms $e : A \otimes A^\circ \to I$ and $n : I \to A^\circ \otimes A$, extraordinary pseudonatural, and so on (for example, see [10]). It also follows that we have a pseudonatural equivalence

$$\mathcal{M}(B \otimes A^\circ, C) \simeq \mathcal{M}(B, C \otimes A).$$

A monoidale $E$ is called *right coclosed* when the composite

$$\hat{p} = \left( E \xrightarrow{n \otimes 1} E^\circ \otimes E \otimes E \xrightarrow{1 \otimes p} E^\circ \otimes E \right)$$

has a left adjoint.

**Example 3.2.** Let $\mathcal{V}$ be a complete cocomplete symmetric monoidal closed category. Let $\mathcal{V}$-Mod denote the autonomous monoidal bicategory whose objects are $\mathcal{V}$-categories and whose hom categories are defined by

$$\mathcal{V}\text{-Mod}(A, B) = [B^{\mathrm{op}} \otimes A, \mathcal{V}]$$

in the $\mathcal{V}$-functor-category notation of [14]. We call the objects of $\mathcal{V}$-Mod$(A, B)$ *modules* from the $\mathcal{V}$-category $A$ to the $\mathcal{V}$-category $B$ (they are also called "bimodules", "profunctors" and "distributors"). Each $\mathcal{V}$-functor $t : A \to B$ determines a module $t_* : A \to B$ defined by $t_*(b, a) = B(b, t(a))$. In fact, $t_*$ has a right adjoint $t^*$ defined by $t^*(a, b) = B(t(a), b)$. Moreover, autonomy of $\mathcal{V}$-Mod is clear by taking $E^\circ = E^{\mathrm{op}}$. Suppose $C$ is a right-closed monoidal $\mathcal{V}$-category; so we have

$$C(x \otimes y, z) \cong C(y, z^x).$$

We can regard the right internal hom as a $\mathcal{V}$-functor $h : C \otimes C^{\mathrm{op}} \to C^{\mathrm{op}}$, $(x, z) \mapsto z^x$. Then our example is that $E = C^{\mathrm{op}}$ becomes a coclosed monoidale in $\mathcal{V}$-Mod with the module $p : E \otimes E \to E$ defined by

$$p(z, x, y) = E(z, x \otimes y) = C(x \otimes y, z).$$

This is because $\hat{p}(x, z, y) = p(z, x, y) \cong E(z^x, y) = h^*(x, z, y)$, so that $\hat{p} = h^*$ has the left adjoint $h_*$.

Suppose $A$ and $E$ are monoidales in $\mathcal{M}$. Dualizing a concept from Section 9 of [12], we say that an opmonoidal morphism $w : A \to E$, with a right adjoint $w^*$ in $\mathcal{M}$, is *strong right coclosed* when $\psi : w\, p \Rightarrow p(w \otimes w)$ has an invertible mate:

$$\psi^r : p(1 \otimes w^*) \Rightarrow w^*\, p(w \otimes 1).$$

## 4. Monoidal comonads

When $\mathcal{M}$ is a monoidal bicategory, so too is Comnd$\mathcal{M}$ where

$$(A, g) \otimes (B, h) = (A \otimes B, g \otimes h).$$

A *monoidal comonad* in $\mathcal{M}$ is a comonad $g = (g, \phi, \phi_0) : E \to E$ in the bicategory Mon$\mathcal{M}$. Equally, it is a monoidale in the monoidal bicategory Comnd$\mathcal{M}$. It consists of

| $E$ | $g$ | $E$ | $\phi : p(g \otimes g) \Rightarrow g\, p$ | $\phi_0 : j \Rightarrow g\, j$ |
|---|---|---|---|---|
| $(p, \phi) : (E \otimes E, g \otimes g) \to (E, g)$ | | $(j, \phi_0) : (I, 1) \to (E, g)$ | | Comnd$\mathcal{M}$ |

$$\mathcal{M} \qquad \text{Chikhladze-Lack-Street} \qquad g = (g, \phi, \phi_0) : E \to E$$

a monoidale $E$, a comonad $g$ on $E$, and 2-cells $\phi : p(g \otimes g) \Rightarrow g\, p$, $\phi_0 : j \Rightarrow g\, j$ such that $(p, \phi) : (E \otimes E, g \otimes g) \to (E, g)$ and $(j, \phi_0) : (I, 1) \to (E, g)$ are morphisms in Comnd$\mathcal{M}$.

If $\mathcal{M}$ admits the Eilenberg-Moore construction for comonads then so too does Mon$\mathcal{M}$; see [15] and [16]. We now generalize a definition due to Bruguières-Lack-Virelizier announced in [1]. All of Proposition 4.2, Theorem 4.3 and Proposition 4.4 below are generalizations to our setting of the results forming part of [1-2].

**Definition 4.1.** A monoidal comonad $g$ on a monoidal $E$ is *right Hopf* when the pasted 2-cell 4.1 is invertible.

$$(4.1)$$

$$
\begin{array}{ccccc}
E \otimes E & \xrightarrow{\;g \otimes 1\;} & E \otimes E & \xrightarrow{\;p\;} & E \\
{\scriptstyle 1 \otimes g}\downarrow & \overset{\delta \otimes 1_g}{\Longrightarrow} & {\scriptstyle g \otimes g}\downarrow & \overset{\phi}{\Longrightarrow} & \downarrow{\scriptstyle g} \\
E \otimes E & \xrightarrow[\;g \otimes 1\;]{} & E \otimes E & \xrightarrow[\;p\;]{} & E
\end{array}
$$

We henceforth assume that the monoidal bicategory $\mathcal{M}$ is right autonomous and admits the Eilenberg-Moore construction for comonads.

**Proposition 4.2.** *A monoidal comonad $g$ on $E$ is right Hopf if and only if the 2-cell 4.2 is invertible.*

$$(4.2)$$

$$
\begin{array}{ccccc}
E^g \otimes E & \xrightarrow{\;u \otimes 1\;} & E \otimes E & \xrightarrow{\;p\;} & E \\
{\scriptstyle 1 \otimes g}\downarrow & \overset{\gamma \otimes 1_g}{\Longrightarrow} & {\scriptstyle g \otimes g}\downarrow & \overset{\phi}{\Longrightarrow} & \downarrow{\scriptstyle g} \\
E^g \otimes E & \xrightarrow[\;u \otimes 1\;]{} & E \otimes E & \xrightarrow[\;p\;]{} & E
\end{array}
$$

*Proof.* Consider the following $3 \times 3$ diagram.

$$(4.3)$$

$$
\begin{array}{ccccc}
u \otimes g & \xrightarrow{\;\gamma \otimes 1_g\;} & g\, u \otimes g & \xrightarrow{\;\phi(u \otimes 1)\;} & g(u \otimes 1) \\
{\scriptstyle \gamma \otimes 1}\downarrow & = & \downarrow{\scriptstyle g\, \gamma \otimes 1} & = & \downarrow{\scriptstyle g(\gamma \otimes 1)} \\
g\, u \otimes g & \xrightarrow{\;\delta\, u \otimes 1_g\;} & g^2\, u \otimes g & \xrightarrow{\;\phi(g\, u \otimes 1)\;} & g(g\, u \otimes 1) \\
{\scriptstyle g\ \ \gamma \otimes 1}\downarrow\ \downarrow{\scriptstyle \delta\, u \otimes 1} & = & {\scriptstyle g^2\, \gamma \otimes 1}\downarrow\ \downarrow{\scriptstyle g\, \delta\, u \otimes 1} & = & {\scriptstyle g(g\, \gamma \otimes 1)}\downarrow\ \downarrow{\scriptstyle g(\delta\, u \otimes 1)} \\
g^2\, u \otimes g & \xrightarrow[\;\delta\, g\, u \otimes 1_g\;]{} & g^3\, u \otimes g & \xrightarrow[\;\phi(g^2\, u \otimes 1)\;]{} & g(g^2\, u \otimes 1)
\end{array}
$$

Since the pseudofunctor $- \otimes E : \mathcal{M} \to \mathcal{M}$ has a left biadjoint $- \otimes E^\circ$, it preserve Eilenberg-Moore constructions (as they are weighted limits). It follows that $u \otimes E : E^g \otimes E \longrightarrow E \otimes E$ is comonadic and then that the columns of 4.3 are

$$u^* \otimes E$$

$$-\otimes E : \mathcal{M} \to \mathcal{M} \qquad\qquad -\otimes E^{\circ}$$

$$u \otimes E : E^g \otimes E \;\longrightarrow\; E \otimes E$$

absolute equalizers. If 4.1 is invertible then the composites in the second and third row of 4.3 are invertible. Therefore the composite in the first row is invertible yielding the invertibility of 4.2. For the converse, note that precomposing 4.2 with $u^* \otimes E$ yields 4.1. $\quad\square$

**Theorem 4.3.** *Suppose $E$ is a right coclosed monoidale in $\mathcal{M}$. A monoidal comonad $g$ on $E$ is right Hopf if and only if the monoidale $E^g$ is right coclosed and $u : E^g \to E$ is strong right coclosed.*

*Proof.* Suppose $g$ is a Hopf monoidal comonad on $E$. Consider the following diagram.

$$(4.4)$$

The pasted 2-cell is invertible by Proposition 4.2 while the top and bottom composites are isomorphic to the composite

$$(4.5) \qquad\qquad E \xrightarrow{\hat{p}} E^{\circ} \otimes E \xrightarrow{u^{\circ} \otimes 1} E^{g\,\circ} \otimes E.$$

Now $\hat{p}$ has a left adjoint since $E$ is right coclosed and $u^{\circ} \otimes 1$ has the left adjoint $u^* \otimes 1$. By Corollary 2.2, the induced morphism $q : E^g \to E^{g\,\circ} \otimes E^g$ on the Eilenberg-Moore constructions of the comonads at the left and right ends of diagram 4.4 has a left adjoint. In fact, $q = \hat{p}$ where this $p$ is that of the monoidale $E^g$. To see this we must see that pasting the triangle containing the universal $g$-coalgebra $\gamma$ to the left of diagram 4.4, conjugated by the isomorphisms of the top and bottom with the morphism 4.5, is equal to $\hat{p}$ pasted on the right with the triangle for the universal $\left(E^{g\,\circ} \otimes g\right)$-coalgebra $E^{g\,\circ} \otimes \gamma$. The morphism part of this calculation is:

$$\left(E^{g\,\circ} \otimes u\right)\hat{p} = \left(E^{g\,\circ} \otimes u\right)\left(E^{g\,\circ} \otimes p\right)\left(n \otimes E^g\right)$$

$$\cong \left(E^{g\,\circ} \otimes p\right)\left(u \otimes u\right)\left(n \otimes E^g\right) \cong \left(E^{g\,\circ} \otimes p\right)\left(u^{\circ} \otimes E \otimes E\right)\left(n \otimes E\right)u$$

$$\cong \left(u^{\circ} \otimes E\right)\left(E^{\circ} \otimes p\right)\left(n \otimes E\right)u \cong \left(u^{\circ} \otimes E\right)\hat{p}\,u\,.$$

For the 2-cell part of the calculation notice that monoidalness of $\gamma$ implies the equality of the composite 2-cells

$$(1 \otimes p)\left(1 \otimes u \otimes u\right)\left(n \otimes 1\right) \xRightarrow{(1 \otimes p)\,(1 \otimes \gamma \otimes \gamma)\,(n \otimes 1)}$$
$$(1 \otimes p)\left(1 \otimes g\,u \otimes g\,u\right)\left(n \otimes 1\right) = (1 \otimes p)\left(1 \otimes g \otimes g\right)\left(1 \otimes u \otimes u\right)\left(n \otimes 1\right)$$
$$\xRightarrow{\psi(1 \otimes u \otimes u)\,(n \otimes 1)} (1 \otimes g)(1 \otimes p)\left(1 \otimes u \otimes u\right)\left(n \otimes 1\right) \cong$$
$$(1 \otimes g)(1 \otimes u)(1 \otimes p)\left(n \otimes 1\right)$$

and

$$(1 \otimes p)\left(1 \otimes u \otimes u\right)\left(n \otimes 1\right) \cong (1 \otimes u)(1 \otimes p)\left(n \otimes 1\right) \xRightarrow{\gamma(1 \otimes p)\,(n \otimes 1)}$$
$$(1 \otimes g\,u)(1 \otimes p)\left(n \otimes 1\right) = (1 \otimes g)(1 \otimes u)(1 \otimes p)\left(n \otimes 1\right)\,.$$

Hence it follows that $E^g$ is right coclosed.

We need to see that $u : E^g \to E$ is strong right coclosed. That is, we need to show $p(1 \otimes u^*) \Rightarrow u^* \, p(u \otimes 1)$ is invertible. As the unit for the adjunction $u \dashv u^*$ is an equalizer, composition with $u$ is conservative (invertibility reflecting). So it suffices to show that $u \, p(1 \otimes u^*) \Rightarrow u \, u^* \, p(u \otimes 1)$. But, using the invertibility of 4.2, we have

$$u \, p(1 \otimes u^*) \cong p(u \otimes u) \, (1 \otimes u^*) \cong p(u \otimes g) \cong g \, p(u \otimes 1) \cong u \, u^* \, p(u \otimes 1).$$

The converse follows from the more general fact expressed in Proposition 4.4.  $\square$

**Proposition 4.4.** *If $w : A \to E$ is a strong right coclosed, strong monoidal morphism with a right adjoint $w^*$ in $\mathcal{M}$ then the generated monoidal comonad $g = w \, w^*$ on $E$ is Hopf.*

*Proof.* We have

$$p(g \otimes 1) \, (1 \otimes g) \cong p(g \otimes g) \cong p(w \otimes w) \, (w^* \otimes w^*)$$

and, using that $w$ is strong monoidal, we continue with

$$\cong w \, p(w^* \otimes w^*) \cong w \, p(1 \otimes w^*) \, (w^* \otimes 1)$$

and, using that $w$ is strong right closed, we continue with

$$\cong w \, w^* \, p(w \otimes 1) \, (w^* \otimes 1) \cong w \, w^* \, p(w \, w^* \otimes 1) \cong g \, p(g \otimes 1).  \quad \square$$

**Example 4.5.** In the setting of Example 3.2, suppose $t = (t, \mu, \eta, \psi, \psi_0)$ is an opmonoidal monad on a monoidal $\mathcal{V}$-category $C$. We may consider $t$ equally to be a monoidal comonad, which we will write as $g = (g, \delta, \varepsilon, \phi, \phi_0)$ to mark the distinction, on $E = C^{\mathrm{op}}$. Then $t$ is a Hopf opmonoidal monad on $C$ in the sense of [1] (actually, they say more briefly that $t$ is a "Hopf monad") if and only if $g_*$ is a Hopf monoidal comonad on $E = C^{\mathrm{op}}$ in $\mathcal{V}$-Mod. To be explicit, the right Hopf condition is, for all objects $x$ and $y$ in $C$, the invertibility of the composite

(4.6) $$\qquad t(t(x) \otimes y) \xrightarrow{\psi_{t(x),y}} t^2(x) \otimes t(y) \xrightarrow{\mu_x \otimes 1_{t(x)}} t(x) \otimes t(y).$$

It should be noted that, if $C^t$ denotes the $\mathcal{V}$-category of Eilenberg-Moore $t$-algebras, then $E^g = (C^t)^{\mathrm{op}}$ in $\mathcal{V}$-Mod. By Proposition 4.2, the invertibility of the arrows 4.6 is equivalent to the invertibility of the arrows

(4.7) $$\qquad t(a \otimes y) \xrightarrow{\psi_{a,y}} t(a) \otimes t(y) \xrightarrow{\alpha \otimes 1_{t(x)}} a \otimes t(y)$$

for all $t$-algebras $(a, \alpha)$ and all objects $y$. If $C$ is right closed, it follows that, for all $t$-algebras $(a, \alpha)$ and $(b, \beta)$, there is a unique $t$-action $\rho : t(b^a) \to b^a$ on $b^a$ such that

$$\left( t(a \otimes b^a) \xrightarrow{t(\mathrm{ev})} t(b) \xrightarrow{\beta} b \right) = \left( t(a \otimes b^a) \xrightarrow{t(\mathrm{ev})} t(a) \otimes t(b^a) \longrightarrow a \otimes b^a \xrightarrow{\mathrm{ev}} b \right).$$

As pointed out in [2], this means that $\mathrm{ev} : a \otimes b^a \longrightarrow b$ becomes a $t$-algebra morphism. Hence $(b^a, \rho)$ is a right internal hom of $(a, \alpha)$ and $(b, \beta)$ in $C^t$, all in accord with our Theorem 4.3.

## 5. Quantum groupoids

Let $\mathcal{V}$ be a braided monoidal category with coreflexive equalizers preserved by the functors $X \otimes - : \mathcal{V} \to \mathcal{V}$. Following [17], we consider the right autonomous monoidal bicategory $\mathrm{Comod}(\mathcal{V})$ whose objects are comonoids $C = (C, \delta)$ (sometimes called "coalgebras") in $\mathcal{V}$ and whose hom $\mathrm{Comod}(\mathcal{V})(C, D)$ is the Eilenberg-Moore category for the comonad $C \otimes - \otimes D$ on $\mathcal{V}$. Eilenberg-Moore coalgebras $M$ for $C \otimes - \otimes D$ are called *comodules from $C$ to $D$* and depicted by $M : C \to D$. The composite $N \circ M = M \diamond_D N$ of comodules $M : C \to D$ and $N : D \to E$ is defined as a coreflexive equalizer

$$M \diamond_D N \longrightarrow M \otimes N \; \underset{\longrightarrow}{\overset{\longrightarrow}{}} \; M \otimes D \otimes N.$$

The monoidal structure on $\mathrm{Comod}(\mathcal{V})$ simply extends the tensor product on $\mathcal{V}$. The right dual $C^\circ$ of the comonoid $C$ has comultiplication defined using the braiding thus:

$$\left( C^\circ \overset{\delta}{\longrightarrow} C^\circ \otimes C^\circ \right) = \left( C \overset{\delta}{\longrightarrow} C \otimes C \overset{c_{C,C}}{\longrightarrow} C \otimes C \right).$$

It follows that the biduality $C \dashv_b C^\circ$ generates a monoidal structure $p = 1 \otimes e \otimes 1$, $j = n$ on $C^\circ \otimes C$ in $\mathrm{Comod}(\mathcal{V})$.

Recall that a *quantum category $A = (C, G)$ in $\mathcal{V}$* is defined in [12] to be a comonoid $C$ in $\mathcal{V}$ equipped with a monoidal comonad $G$ on the monoidale $C^\circ \otimes C$ in $\mathrm{Comod}(\mathcal{V})$. The authors go on to add conditions, involving star-autonomy, in order for a quantum category to be a quantum groupoid. Now we shall suggest a more general definition which is in the spirit of [10], [18] and [19]. A virtue is that it is a property, rather than extra structure, that a quantum category may have.

**Definition 5.1.** A *quantum groupoid in $\mathcal{V}$* is a quantum category $A = (C, G)$ where the monoidal comonad $G$ is Hopf.

## References

[1] Alain Bruguières, "Hopf monads on monoidal categories, and modularity of the centre" *Australian Category Seminar* (Macquarie University, 10 December 2008); http://www.maths.usyd.edu.au/u/AusCat/abstracts/081210ab.html

[2] A. Bruguières, S. Lack and A. Virelizier (2009). Hopf monads on monoidal categories. in preparation

[3] S. Eilenberg and G. M. Kelly, Closed categories, *Proc. Conf. Categorical Algebra (La Jolla, Calif., 1965)*, Springer-Verlag, New York, 1966, 421–562. MR **37**:1432

[4] A. Bruguières and A. Virelizier, Hopf monads, Advances in Mathematics **215** (2007), 679–733.

[5] C. Pastro and R. Street, Closed categories, star-autonomy, and monoidal comonads, Journal of Algebra **321** (2009), 3494–3520.

[6] J. Bénabou, Introduction to bicategories, Lecture Notes in Mathematics **47** (Springer-Verlag, Berlin, 1967), 1–77.

[7] R. Street, The formal theory of monads, Journal of Pure and Applied Algebra **2** (1972), 149–168.

[8] G. M. Kelly, Doctrinal adjunction, Category Seminar (Proc. Sem., Sydney, 1972/1973) Lecture Notes in Math., **420** (Springer, Berlin,1974), 257–280.

[9] R. Street, Two constructions on lax functors, Cahiers topologie et géométrie différentielle **13** (1972), 217–264.

[10] B. Day and R. Street, Monoidal bicategories and Hopf algebroids, Advances in Mathematics **129** (1997), 99–157.

[11] R. Gordon, A. J. Power and R. Street, *Coherence for tricategories*, Memoirs of the American Math. Society, Volume117, Number 558, 1995.

[12] B. Day and R. Street, Quantum categories, star autonomy, and quantum groupoids, Fields Institute Communications (American Math. Soc.) **43** (2004), 187–226.

[13] I. Lopez-Franco, R. Street and R. J. Wood, Duals invert, Applied Categorical Structures **?** (accepted 8 September 2009; published online: 26 September 2009 ), x–x+40. http://www.springerlink.com/content/u517174h46wm2313/

[14] G. M. Kelly, *Basic concepts of enriched category theory*, London Mathematical Society Lecture Note Series, 64. Cambridge University Press, Cambridge, 1982. MR **84e**:18001 http://www.tac.mta.ca/tac/reprints/articles/10/tr10.pdf

[15] I. Moerdijk, Monads on tensor categories, Journal of Pure Appl. Algebra **168** (2002), 189–208.

[16] P. McCrudden, Opmonoidal monads, Theory and Applications of Categories **10** (2002), 469–485. http://www.tac.mta.ca/tac/volumes/10/19/10-19.pdf

[17] B. Day, P. McCrudden and R. Street, Dualizations and antipodes, Applied Categorical Structures **11** (2003), 229–260.

[18] R. Street, Fusion operators and cocycloids in monoidal categories, Applied Categorical Structures **6** (1998), 177–191.

[19] P. Schauenburg, Duals and doubles of quantum groupoids ($\times$_R-Hopf algebras), Contemporary Mathematics ( Amer. Math. Soc., Providence, RI) **267** (2000), 273–299.